\documentclass[12pt]{article}

\usepackage[geometry]{ifsym}
\usepackage{verbatim}
\usepackage{fancyhdr}
\usepackage{graphicx}
\usepackage{mathrsfs}
\usepackage{amssymb}
\usepackage{pst-poly}

\newtheorem{thm}{Theorem}[section]

\newenvironment{pf}[1][Proof]{\noindent\textbf{#1.} }{\hfill\rule{1mm}{2mm}}

\makeatletter \@addtoreset{equation}{section} \makeatother

\begin{document}

\title{Complexity of Bondage and Reinforcement\thanks{The work was
supported by NNSF of China (No.10671191).}}
\author{Fu-Tao Hu, \quad
Jun-Ming Xu\footnote{\ Correspondence to: J.-M. Xu; e-mail: xujm@ustc.edu.cn}\\
{\small Department of Mathematics} \\
{\small University of Science and Technology of China}\\
{\small Hefei, 230026, China} }

\date{}
\maketitle

\begin{abstract}

\setlength{\baselineskip}{24pt}

Let $G=(V,E)$ be a graph. A subset $D\subseteq V$ is a dominating
set if every vertex not in $D$ is adjacent to a vertex in $D$. A
dominating set $D$ is called a total dominating set if every vertex
in $D$ is adjacent to a vertex in $D$. The domination (resp. total
domination) number of $G$ is the smallest cardinality of a
dominating (resp. total dominating) set of $G$. The bondage (resp.
total bondage) number of a nonempty graph $G$ is the smallest number
of edges whose removal from $G$ results in a graph with larger
domination (resp. total domination) number of $G$. The reinforcement
number of $G$ is the smallest number of edges whose addition to $G$
results in a graph with smaller domination number. This paper shows
that the decision problems for bondage, total bondage and
reinforcement are all NP-hard.


\end{abstract}

\newpage

\setlength{\baselineskip}{24pt}

\section{Introduction}

In this paper, we follow Xu~\cite{xu03} for graph-theoretical
terminology and notation. A graph $G=(V,E)$ always means a finite,
undirected and simple graph, where $V=V(G)$ is the vertex-set and
$E=E(G)$ is the edge-set of $G$.

A subset $D\subseteq V$ is a {\it dominating set} of $G$ if every
vertex not in $D$ is adjacent to a vertex in $D$. The {\it
domination number} of $G$, denoted by $\gamma(G)$, is the minimum
cardinality of a dominating set of $G$. A dominating set $D$ is
called a {\it $\gamma$-set} of $G$ if $|D|=\gamma(G)$. The {\it
bondage number} of $G$, denoted by $b(G)$, is the minimum number of
edges whose removal from $G$ results in a graph with larger
domination number of $G$. The {\it reinforcement number} of $G$,
denoted by $r(G)$, is the smallest number of edges whose addition to
$G$ results in a graph with smaller domination number of $G$.
Domination is a classical concept in graph theory. The bondage
number and the reinforcement number were introduced by Fink et
at.~\cite{fjkr90} and Kok, Mynhardt~\cite{km90}, respectively, in
1990. The reinforcement number for digraphs has been studies by
Huang, Wang and Xu~\cite{hwx09}. Domination as well as related
topics is now well studied in graph theory. The literature on these
subjects have been surveyed and detailed in the two excellent
domination books by Haynes, Hedetniemi, and Slater~\cite{hh197,
hh297}.

Theory of domination has been applied in many research fields. For
different applications, many variations of dominations were proposed
in the research literature by adding some restricted conditions to
dominating sets, for example, the total domination and the
restrained domination.

A dominating set $D$ is called a {\it total dominating set} if every
vertex in $D$ is adjacent to another vertex in $D$. The {\it total
domination number}, denoted by $\gamma_t(G)$, of $G$ is the minimum
cardinality of a total dominating set of $G$. Use the symbol $D_t$
to denote a total dominating set. A total dominating set $D_t$ is
called a {\it $\gamma_t$-set} of $G$ if $|D_t|=\gamma_t(G)$. The
{\it total bondage number} of $G$, denoted by $b_t(G)$, is the
minimum number of edges whose removal from $G$ results in a graph
with larger total domination number of $G$. The total domination was
introduced by Cockayne et al.~\cite{cdh80}. Total domination in
graphs has been extensively studied in the literature. A survey of
selected recent results on total domination in Henning~\cite{h09}.
The total bondage number of a graph was first studied by Kulli and
Patwari~\cite{kp91} and further studied by Sridharan, Elias,
Subramanian~\cite{ses07}, Huang and Xu~\cite{hx07a}.

Analogously, a dominating set $D$ is called a {\it restrained
dominating set} if every vertex not in $D$ is adjacent to another
vertex not in $D$. The {\it restrained domination number}, denoted
by $\gamma_r(G)$, of $G$ is the minimum cardinality of a total
dominating set of $G$. The {\it restrained bondage number} of $G$,
denoted by $b_r(G)$, is the minimum number of edges whose removal
from $G$ results in a graph with larger restrained domination number
of $G$. The restrained domination was introduced by Telle and
Proskurowski~\cite{tp97}, and the restrained bondage number was
defined by Hattingh and Plummer~\cite{hp08}.

Whys that a graph-theoretical parameter is proposed at once is to
determine the exact value of this parameter for all graphs. However,
the problem determining domination for general graphs has been
proved to be NP-complete (see GT2 in Appendix in Garey and
Johnson~\cite{gj79}); the problems determining total domination and
restrained domination for general graphs have been also proved to be
NP-complete by Laskar et al.~\cite{lphh84}, and by Domke et
at.~\cite{dh99}, respectively.

As regards the bondage problem, Hattingh et al.~\cite{hp08} showed
that the restrained bondage problem is NP-complete even for
bipartite graphs. For the general bondage problem, from the
algorithmic point of view, Hartnell et at.~\cite{hjvw98} designed a
linear time algorithm to compute the bondage number of a tree.
However, the complexity of this problem is still unknown for other
classes of graphs.

In this paper, we will show that the decision problems for bondage,
total bondage and reinforcement are all NP-hard. Their proofs
are Section 3, Section 4 and Section 5 in this paper, respectively.

\section{$3$-satisfiability problem}

Following Garey and Johnson's techniques for proving
NP-hardness~\cite{gj79}, we prove our results by describing a
polynomial transformation from the known NP-complete problem:
$3$-satisfiability problem. To state the $3$-satisfiability problem,
we, in this section, recall some terms we will use in describing it.

Let $U$ be a set of Boolean variables. A {\it truth assignment} for
$U$ is a mapping $t: U\to\{T,F\}$. If $t(u)=T$, then $u$ is said to
be ``\,true" under $t$; if If $t(u)=F$, then $u$ is said to
be``\,false" under $t$. If $u$ is a variable in $U$, then $u$ and
$\bar{u}$ are {\it literals} over $U$. The literal $u$ is true under
$t$ if and only if the variable $u$ is true under $t$; the literal
$\bar{u}$ is true if and only if the variable $u$ is false.

A {\it clause} over $U$ is a set of literals over $U$. It represents
the disjunction of these literals and is {\it satisfied} by a truth
assignment if and only if at least one of its members is true under
that assignment. A collection $\mathscr C$ of clauses over $U$ is
{\it satisfiable} if and only if there exists some truth assignment
for $U$ that simultaneously satisfies all the clauses in $\mathscr
C$. Such a truth assignment is called a {\it satisfying truth
assignment} for $\mathscr C$. The $3$-satisfiability problem is
specified as follows.

\begin{center}
\begin{minipage}{130mm}
\setlength{\baselineskip}{24pt}

\vskip6pt\noindent {\bf $3$-satisfiability problem}:

\noindent {\bf Instance:}\ {\it A
collection $\mathscr{C}=\{C_1,C_2,\ldots,C_m\}$ of clauses over a
finite set $U$ of variables such that $|C_j| =3$ for $j=1,
2,\ldots,m$.}

\noindent {\bf Question:}\ {\it Is there a truth assignment for $U$
that satisfies all the clauses in $\mathscr{C}$?}

\end{minipage}
\end{center}

\begin{thm} \textnormal{(Theorem 3.1 in~\cite{gj79})}
The $3$-satisfiability problem is NP-complete.
\end{thm}




\section{NP-hardness of bondage}

In this section, we will show that the problem determining the
bondage numbers of general graphs is NP-hard. We first state the
problem as the following decision problem.

\begin{center}
\begin{minipage}{130mm}
\setlength{\baselineskip}{24pt}

\vskip6pt\noindent {\bf Bondage problem:}

\noindent {\bf Instance:}\ {\it A nonempty graph $G$ and a positive
integer $k$.}

\noindent {\bf Question:}\ {\it Is $b(G)\le k$?}

\end{minipage}
\end{center}

\vskip6pt\begin{thm} The bondage problem is NP-hard.
\end{thm}

\begin{pf}
We show the NP-hardness of the bondage problem by
transforming the $3$-satisfiability problem to it in polynomial
time.

Let $U=\{u_1,u_2,\ldots,u_n\}$ and $\mathscr{C}=\{C_1,C_2,
\ldots,C_m\}$ be an arbitrary instance of the $3$-satisfiability
problem. We will construct a graph $G$ and a positive integer $k$
such that $\mathscr{C}$ is satisfiable if and only if $b(G)\leq k$.
Such a graph $G$ can be constructed as follows.

For each $i=1,2,\ldots,n$, corresponding to the variable $u_i\in U$,
associate a triangle $T_i$ with vertex-set $\{u_i,\bar{u_i},v_i\}$.
For each $j=1,2,\ldots,m$, corresponding to the clause
$C_j=\{x_j,y_j,z_j\}\in \mathscr{C}$, associate a single vertex
$c_j$ and add edge-set $E_j=\{c_jx_j, c_jy_j, c_jz_j\}$. Finally,
add a path $P=s_1s_2s_3$, join $s_1$ and $s_3$ to each vertex $c_j$
with $1\le j\le m$ and set $k=1$.

Figure~\ref{f1} shows an example of the graph obtained when
$U=\{u_1,u_2,u_3,u_4\}$ and $\mathscr{C}=\{C_1,C_2,C_3\}$, where
$C_1=\{u_1,u_2,\bar{u_3}\}, C_2=\{\bar{u_1},u_2,u_4\},
C_3=\{\bar{u_2},u_3,u_4\}$.

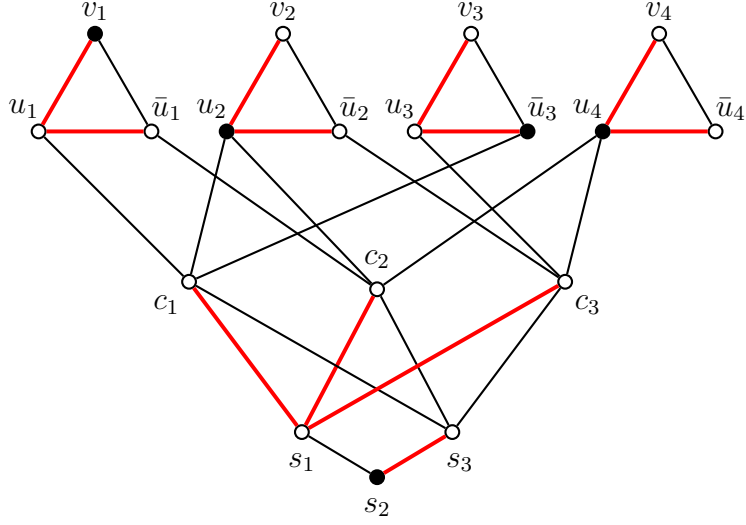
\begin{figure}[ht]
\begin{center}
\begin{pspicture}(-5,-1.1)(5,6.7)

\cnode*(0,-.6){3pt}{s2}\rput(0,-1){$s_2$}
\cnode(-1,0){3pt}{s1}\rput(-1,-.4){$s_1$}
\cnode(1,0){3pt}{s3}\rput(1.1,-.4){$s_3$} \ncline{s2}{s1}
\ncline[linecolor=red,linewidth=1.5pt]{s2}{s3}

\cnode(0,1.9){3pt}{c2}\rput(0,2.3){$c_2$}
\ncline[linecolor=red,linewidth=1.5pt]{c2}{s1} \ncline{c2}{s3}
\cnode(-2.5,2){3pt}{c1}\rput(-2.8,1.7){$c_1$}
\ncline[linecolor=red,linewidth=1.5pt]{c1}{s1} \ncline{c1}{s3}
\cnode(2.5,2){3pt}{c3}\rput(2.8,1.7){$c_3$}
\ncline[linecolor=red,linewidth=1.5pt]{c3}{s1} \ncline{c3}{s3}

\cnode(-4.5,4){3pt}{u1}\rput(-4.7,4.3){$u_1$}
\cnode(-3,4){3pt}{u1'}\rput(-2.8,4.3){$\bar{u}_1$}
\ncline[linecolor=red,linewidth=1.5pt]{u1}{u1'}
\cnode*(-2,4){3pt}{u2}\rput(-2.2,4.3){$u_2$}
\cnode(-0.5,4){3pt}{u2'}\rput(-0.3,4.3){$\bar{u}_2$}
\ncline[linecolor=red,linewidth=1.5pt]{u2}{u2'}
\cnode(0.5,4){3pt}{u3}\rput(0.3,4.3){$u_3$}
\cnode*(2,4){3pt}{u3'}\rput(2.2,4.3){$\bar{u}_3$}
\ncline[linecolor=red,linewidth=1.5pt]{u3}{u3'}
\cnode*(3,4){3pt}{u4}\rput(2.8,4.3){$u_4$}
\cnode(4.5,4){3pt}{u4'}\rput(4.7,4.3){$\bar{u}_4$}
\ncline[linecolor=red,linewidth=1.5pt]{u4}{u4'}
\cnode*(-3.75,5.3){3pt}{v1}\rput(-3.75,5.6){$v_1$}
\ncline[linecolor=red,linewidth=1.5pt]{v1}{u1} \ncline{v1}{u1'}
\cnode(-1.25,5.3){3pt}{v2}\rput(-1.25,5.6){$v_2$}
\ncline[linecolor=red,linewidth=1.5pt]{v2}{u2} \ncline{v2}{u2'}
\cnode(1.25,5.3){3pt}{v3}\rput(1.25,5.6){$v_3$}
\ncline[linecolor=red,linewidth=1.5pt]{v3}{u3} \ncline{v3}{u3'}
\cnode(3.75,5.3){3pt}{v4}\rput(3.75,5.6){$v_4$}
\ncline[linecolor=red,linewidth=1.5pt]{v4}{u4} \ncline{v4}{u4'}

\ncline{c1}{u1} \ncline{c1}{u2} \ncline{c1}{u3'}
\ncline{c2}{u1'} \ncline{c2}{u2} \ncline{c2}{u4}
\ncline{c3}{u2'} \ncline{c3}{u3} \ncline{c3}{u4}
\end{pspicture}
\caption{\label{f1}\footnotesize An instance of the bondage problem
resulting from an instance of the $3$-satisfiability problem, in
which $U=\{u_1,u_2,u_3,u_4\}$ and
$\mathscr{C}=\{\{u_1,u_2,\bar{u_3}\},\{\bar{u_1},u_2,u_4\},
\{\bar{u_2},u_3,u_4\}\}$. Here $k=1$ and $\gamma=5$, where the set
of bold points is a $\gamma$-set.}
\end{center}
\end{figure}

To prove that this is indeed a transformation, we must show that
$b(G)=1$ if and only if there is a truth assignment for $U$ that
satisfies all the clauses in $\mathscr{C}$.
This aim can be obtained by proving the following four claims.

\begin{description}

\item [Claim 3.1]
{\it $\gamma(G)\geq n+1$. Moreover, if $\gamma(G)=n+1$, then for any
$\gamma$-set $D$ in $G$, $D\cap V(P)=\{s_2\}$ and $|D\cap V(T_i)|=1$
for each $i=1,2,\ldots,n$, while $c_j\notin D$ for each
$j=1,2,\ldots,m$.}

\begin{pf}
Let $D$ be a $\gamma$-set of $G$. By the construction of $G$, the
vertex $s_2$ can be dominated only by vertices in $P$, which implies
$|D\cap V(P)|\geq 1$; for each $i=1,2,\ldots,n$, the vertex $v_i$
can be dominated only by vertices in $T_i$, which implies $|D\cap
V(T_i)|\geq 1$. It follows that $\gamma(G)=|D|\geq n+1$.

Suppose that $\gamma(G)=n+1$. Then $|D\cap V(P)|=1$ and $|D\cap
V(T_i)|=1$ for each $i=1,2,\ldots,n$. Consequently, $c_j\notin D$
for each $j=1,2,\ldots,m$. If $s_1\in D$, then $|D\cap V(P)|=1$
implies that $D\cap V(P)=\{s_1\}$, and so $s_3$ could not be
dominated by $D$, a contradiction. Hence $s_1\notin D$. Similarly
$s_3\notin D$ and, thus, $D\cap V(P)=\{s_2\}$ since $|D\cap
V(P)|=1$.
\end{pf}

\item [Claim 3.2]
{\it $\gamma(G)= n+1$ if and only if $\mathscr{C}$ is satisfiable.}

\begin{pf}
Suppose that $\gamma(G)=n+1$ and let $D$ be a $\gamma$-set of $G$.
By Claim 3.1, for each $i=1,2,\ldots,n$, $|D\cap V(T_i)|=1$, it
follows that $D\cap V(T_i)=\{u_i\}$ or $D\cap V(T_i)=\{\bar{u_i}\}$
or $D\cap V(T_i )=\{v_i\}$. Define a mapping $t: U\to \{T,F\}$ by
 \begin{equation}\label{e3.1}
 t(u_i)=\left\{
 \begin{array}{l}
 T \ \ {\rm if}\ u_i\in D \ {\rm or} \ v_i\in D, \\
 F \ \ {\rm if}\ \bar {u_i}\in D,
\end{array}
 \right.
 \ i=1,2,\ldots,n.
 \end{equation}

We will show that $t$ is a satisfying truth assignment for
$\mathscr{C}$. It is sufficient to show that every clause in
$\mathscr{C}$ is satisfied by $t$. To this end, we arbitrarily
choose a clause $C_j\in\mathscr{C}$ with $1\le j\le m$.
Since the corresponding vertex $c_j$ in $G$ is adjacent to neither
$s_2$ nor $v_i$ for any $i$ with $1\le i\le n$, there exists some
$i$ with $1\le i\le n$ such that $c_j$ is dominated by $u_i\in D$ or
$\bar{u}_i\in D$. Suppose that $c_j$ is dominated by $u_i\in D$.
Since $u_i$ is adjacent to $c_j$ in $G$, the literal $u_i$ is in the
clause $C_j$ by the construction of $G$. Since $u_i\in D$, it
follows that $t(u_i)=T$ by (\ref{e3.1}), which implies that the
clause $C_j$ is satisfied by $t$. Suppose that $c_j$ is dominated by
$\bar{u}_i\in D$. Since $\bar{u}_i$ is adjacent to $c_j$ in $G$, the
literal $\bar{u}_i$ is in the clause $C_j$. Since $\bar{u}_i\in D$,
it follows that $t(u_i)=F$ by (\ref{e3.1}). Thus, $t$ assigns
$\bar{u}_i$ the truth value $T$, that is, $t$ satisfies the clause
$C_j$. By the arbitrariness of $j$ with $1\le j\le m$, we show that
$t$ satisfies all the clauses in $\mathscr{C}$, that is,
$\mathscr{C}$ is satisfiable.

Conversely, suppose that $\mathscr{C}$ is satisfiable, and let $t:
U\to \{T,F\}$ be a satisfying truth assignment for $\mathscr{C}$.
Construct a subset $D'\subseteq V(G)$ as follows. If $t(u_i)=T$,
then put the vertex $u_i$ in $D'$; if $t(u_i)=F$, then put the
vertex $\bar{u_i}$ in $D'$. Clearly, $|D'|=n$. Since $t$ is a
satisfying truth assignment for $\mathscr{C}$, for each
$j=1,2,\ldots,m$, at least one of literals in $C_j$ is true under
the assignment $t$. It follows that the corresponding vertex $c_j$
in $G$ is adjacent to at least one vertex in $D'$ since $c_j$ is
adjacent to each literal in $C_j$ by the construction of $G$. Thus
$D'\cup \{s_2\}$ is a dominating set of $G$, and so $\gamma(G)\leq
|D'\cup \{s_2\}|=n+1$. By Claim 3.1, $\gamma(G)\geq n+1$, and so
$\gamma(G)=n+1$.
\end{pf}

\item [Claim 3.3]
{\it $\gamma(G-e)\leq n+2$ for any $e\in E(G)$.}

\begin{pf}
Let $E_1=\{s_2s_3,s_1c_j,u_i\bar{u_i},u_iv_i,: i=1,2,\ldots,n;
j=1,2,\ldots,m\}$ (induced by heavy edges in Figure~\ref{f1}) and
let $E_2=E(G)\setminus E_1$. Assume $e\in E_2$. Let
$D'=\{u_1,u_2,\ldots,u_n,s_1, s_2\}$. Clearly, $D'$ is a dominating
set of $G-e$ since every vertex not in $D'$ is incident with some
vertex in $D'$ via an edge in $E_1$. Hence, $\gamma(G-e)\leq
|D'|=n+2$. Now assume $e\in E_1$. Let $D''=\{u_1,u_2,\ldots,u_n,s_2,
s_3\}$. If $e$ is either $s_2s_3$ or incident with the vertex $s_1$,
then $D''$ is a dominating set of $G-e$, clearly. If $e$ is either
$u_i\bar{u_i}$ or $u_iv_i$ for some $i$ ($1\le i\le n$), then we use
the vertex either $v_i$ or $\bar{u_i}$ instead of $u_i$ in $D''$ to
obtain $D'''$; and hence $D'''$ is a dominating set of $G-e$. These
facts imply that $\gamma(G-e)\leq n+2$.
\end{pf}

\item [Claim 3.4]
{\it $\gamma(G)=n+1$ if and only if $b(G)=1$.}

\begin{pf}
Assume $\gamma(G)=n+1$ and consider the edge $e=s_1s_2$. Suppose
$\gamma(G)=\gamma(G-e)$. Let $D'$ be a $\gamma$-set in $G-e$. It is
clear that $D'$ is also a $\gamma$-set of $G$. By Claim 3.1 we have
$c_j\notin D'$ for each $j=1,2,\ldots,m$ and $D'\cap V(P)=\{s_2\}$.
But then $s_1$ is not dominated by $D'$, a contradiction. Hence,
$\gamma(G)<\gamma(G-e)$, and so $b(G)=1$.

Now, assume $b(G)=1$. By Claim 3.1, we have that $\gamma(G)\geq
n+1$. Let $e'$ be an edge such that $\gamma(G)<\gamma(G-e')$. By
Claim 3.3, we have that $\gamma(G-e')\leq n+2$. Thus, $n + 1\leq
\gamma(G)< \gamma(G-e')\leq n+2$, which yields $\gamma(G)=n+1$.
\end{pf}

\end{description}

By Claim 3.2 and Claim 3.4, we prove that $b(G)=1$ if and only if
there is a truth assignment for $U$ that satisfies all the clauses
in $\mathscr{C}$. Since the construction of the bondage instance is
straightforward from a $3$-satisfiability instance, the size of the
bondage instance is bounded above by a polynomial function of the
size of $3$-satisfiability instance. It follows that this is a
polynomial transformation.

The theorem follows.
\end{pf}

\section{NP-hardness of total bondage}

In this section, we will show that the problem determining the total
bondage numbers of general graphs is NP-hard. We first state it
as the following decision problem.

\begin{center}
\begin{minipage}{130mm}
\setlength{\baselineskip}{24pt}

\vskip6pt\noindent {\bf Total bondage problem:}

\noindent {\bf Instance:}\ {\it A nonempty graph $G$ and a positive
integer $k$.}

\noindent {\bf Question:}\ {\it Is $b_t(G)\le k$?}

\end{minipage}
\end{center}

\vskip6pt\begin{thm} The total bondage problem is NP-hard.
\end{thm}

\begin{pf}
We show the
NP-hardness of the total bondage problem by reducing the
$3$-satisfiability problem to it in polynomial time.

Let $U=\{u_1,u_2,\ldots,u_n\}$ and $\mathscr{C}=\{C_1,C_2,
\ldots,C_m\}$ be an arbitrary instance of the $3$-satisfiability
problem. We will construct a graph $G$ and an integer $k$ such that
$\mathscr{C}$ is satisfiable if and only if $b_t(G)\leq k$. Such a
graph $G$ can be constructed as follows.

For each $i=1,2,\ldots,n$, corresponding to the variable $u_i\in U$,
associate a graph $H_i$ with vertex-set
$V(H_i)=\{u_i,\bar{u_i},v_i,v_i'\}$ and edge-set
$E(H_i)=\{v_iu_i,u_i\bar{u_i},\bar{u_i}v_i,v_iv_i'\}$. For each
$j=1,2,\ldots,m$, corresponding to the clause
$C_j=\{x_j,y_j,z_j\}\in \mathscr{C}$, associate a single vertex
$c_j$ and add edge-set $E_j=\{c_jx_j, c_jy_j,c_jz_j\}$, $1\le j\le
m$. Finally, add a graph $H$ with vertex-set
$V(H)=\{s_1,s_2,s_3,s_4,s_5\}$ and edge-set
$E(H)=\{s_1s_2,s_1s_4,s_2s_3,s_2s_4,s_4s_5\}$, join $s_1$ and $s_3$
to each vertex $c_j$, $1\le j\le m$ and set $k=1$.

Figure~\ref{f2} shows an example of the graph obtained when
$U=\{u_1,u_2,u_3,u_4\}$ and $\mathscr{C}=\{C_1,C_2,C_3\}$, where
$C_1=\{u_1,u_2,\bar{u_3}\}, C_2=\{\bar{u_1},u_2,u_4\}$ and $C_3=
\{\bar{u_2},u_3,u_4\}$.

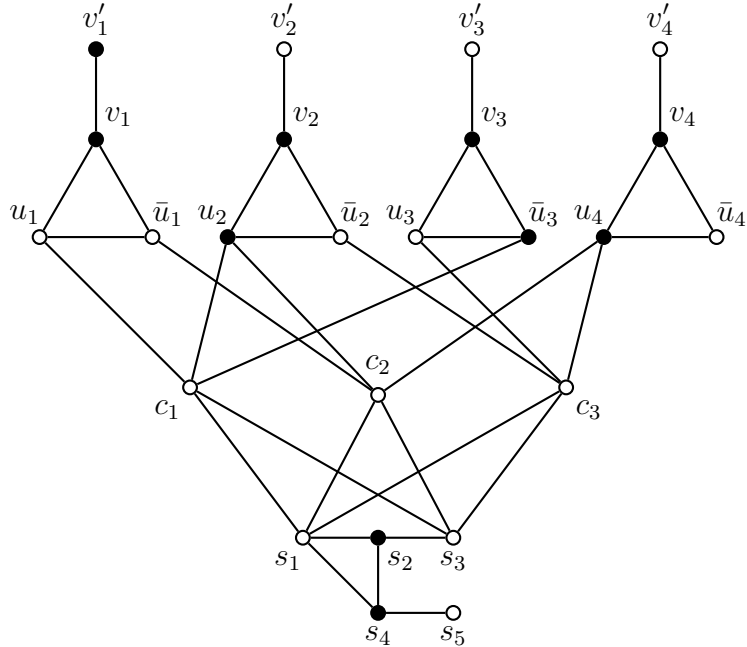
\begin{figure}[ht]
\begin{center}
\begin{pspicture}(-5,-1.5)(5,7.5)

\cnode*(0,0){3pt}{s2}\rput(.3,-.3){$s_2$}
\cnode(-1,0){3pt}{s1}\rput(-1.2,-.3){$s_1$}
\cnode(1,0){3pt}{s3}\rput(1,-.3){$s_3$} \ncline{s2}{s1} \ncline{s2}{s3}
\cnode*(0,-1){3pt}{s4}\rput(0,-1.3){$s_4$} \ncline{s4}{s1} \ncline{s4}{s2}
\cnode(1,-1){3pt}{s5}\rput(1,-1.3){$s_5$} \ncline{s5}{s4}

\cnode(0,1.9){3pt}{c2}\rput(0,2.3){$c_2$} \ncline{c2}{s1}
\ncline{c2}{s3} \cnode(-2.5,2){3pt}{c1}\rput(-2.8,1.7){$c_1$}
\ncline{c1}{s1} \ncline{c1}{s3}
\cnode(2.5,2){3pt}{c3}\rput(2.8,1.7){$c_3$} \ncline{c3}{s1}
\ncline{c3}{s3}

\cnode(-4.5,4){3pt}{u1}\rput(-4.7,4.3){$u_1$}
\cnode(-3,4){3pt}{u1'}\rput(-2.8,4.3){$\bar{u}_1$} \ncline{u1}{u1'}
\cnode*(-2,4){3pt}{u2}\rput(-2.2,4.3){$u_2$}
\cnode(-0.5,4){3pt}{u2'}\rput(-0.3,4.3){$\bar{u}_2$}
\ncline{u2}{u2'} \cnode(0.5,4){3pt}{u3}\rput(0.3,4.3){$u_3$}
\cnode*(2,4){3pt}{u3'}\rput(2.2,4.3){$\bar{u}_3$} \ncline{u3}{u3'}
\cnode*(3,4){3pt}{u4}\rput(2.8,4.3){$u_4$}
\cnode(4.5,4){3pt}{u4'}\rput(4.7,4.3){$\bar{u}_4$} \ncline{u4}{u4'}
\cnode*(-3.75,5.3){3pt}{v1}\rput(-3.45,5.6){$v_1$} \ncline{v1}{u1}
\ncline{v1}{u1'} \cnode*(-1.25,5.3){3pt}{v2}\rput(-.95,5.6){$v_2$}
\ncline{v2}{u2} \ncline{v2}{u2'}
\cnode*(1.25,5.3){3pt}{v3}\rput(1.55,5.6){$v_3$} \ncline{v3}{u3}
\ncline{v3}{u3'} \cnode*(3.75,5.3){3pt}{v4}\rput(4.05,5.6){$v_4$}
\ncline{v4}{u4} \ncline{v4}{u4'}
\cnode*(-3.75,6.5){3pt}{v1'}\rput(-3.75,6.9){$v_1'$}
\ncline{v1}{v1'} \cnode(-1.25,6.5){3pt}{v2'}\rput(-1.25,6.9){$v_2'$}
\ncline{v2}{v2'} \cnode(1.25,6.5){3pt}{v3'}\rput(1.25,6.9){$v_3'$}
\ncline{v3}{v3'} \cnode(3.75,6.5){3pt}{v4'}\rput(3.75,6.9){$v_4'$}
\ncline{v4}{v4'}

\ncline{c1}{u1} \ncline{c1}{u2} \ncline{c1}{u3'}
\ncline{c2}{u1'} \ncline{c2}{u2} \ncline{c2}{u4}
\ncline{c3}{u2'} \ncline{c3}{u3} \ncline{c3}{u4}
\end{pspicture}
\caption{\label{f2}\footnotesize An instance of the total bondage
problem resulting from an instance of the $3$-satisfiability
problem, in which $U=\{u_1,u_2,u_3,u_4\}$ and
$\mathscr{C}=\{\{u_1,u_2,\bar{u}_3\},\{\bar{u}_1,u_2,u_4\},
\{\bar{u}_2,u_3,u_4\}\}$. Here $k=1$ and $\gamma_t=10$, where the
set of bold points is a $\gamma_t$-set.}
\end{center}
\end{figure}

It is easy to see that the construction can be accomplished in
polynomial time. All that remains to be shown is that $\mathscr{C}$
is satisfiable if and only if $b_t(G)=1$. This aim can be obtained
by proving the following four claims.

\begin{description}

\item [Claim 4.1]
{\it $\gamma_t(G)\geq 2n+2$. For any $\gamma_t$-set $D_t$ of $G$,
$s_4\in D_t$ and $v_i\in D_t$ for each $i=1,2,\ldots,n$. Moreover,
if $\gamma_t(G)=2n+2$, then $D_t\cap V(H)=\{s_2,s_4\}$ and $|D_t\cap
V(H_i)|=2$ for each $i=1,2,\ldots,n$, while $c_j\notin D_t$ for each
$j=1,2,\dots,m$.}

\begin{pf}
Let $D_t$ be a $\gamma_t$-set of $G$. By the construction of $G$, it
is clear that $v_i$ is certainly in $D_t$ to dominate $v_i'$, and
$v_i$ can be dominated only by another vertex in $H_i$. It follows
that $v_i\in D_t$ and $|D_t\cap V(H_i)|\geq 2$ for each
$i=1,2,\ldots,n$. It is also clear that $s_4$ is certainly in $D_t$
to dominate $s_5$, and $s_4$ can be dominated only by another vertex
in $H$. This fact implies that $s_4\in D_t$ and $|D_t\cap V(H)|\geq
2$. Thus, $\gamma_t(G)=|D_t|\geq 2n+2$.

Suppose that $\gamma_t(G)=2n+2$. Then $|D_t\cap V(H_i)|=2$ for each
$i=1,2,\ldots,n$, and $|D_t\cap V(H)|=2$. Consequently, $c_j\notin
D_t$ for each $j=1,2,\ldots,m$. As a result, $s_3$ can be dominated
only by the vertex $s_2$ in $S$, that is, $s_2\in D_t$. Noting
$s_4\in D_t$ and $|D_t\cap V(H)|=2$, we have $D_t\cap
V(H)=\{s_2,s_4\}$.
\end{pf}

\item [Claim 4.2]
{\it $\gamma_t(G)=2n+2$ if and only if $\mathscr{C}$ is
satisfiable.}

\begin{pf}
Suppose that $\gamma_t(G)=2n+2$ and let $D_t$ be a $\gamma_t$-set of
$G$. By Claim 4.1, $D_t\cap V(H)=\{s_2,s_4\}$ and for each
$i=1,2,\dots,n$, $|D_t\cap V(H_i)|=2$, it follows that $D_t\cap
V(H_i)=\{u_i,v_i\}$ or $\{\bar{u_i},v_i\}$ or $\{v_i,v_i'\}$. Define
a mapping $t: U\to \{T,F\}$ by
 \begin{equation}\label{e4.1}
 t(u_i)=\left\{
\begin{array}{ll}
 T \ & {\rm if}\ u_i\in D_t \ {\rm or} \ v_i'\in D_t, \\
 F \ & {\rm if}\ \bar {u_i}\in D_t,
\end{array}
 \right.
 \ i=1,2,\ldots, n.
 \end{equation}

We will show that $t$ is a satisfying truth assignment for
$\mathscr{C}$. It is sufficient to show that $t$ satisfies every
clause in $\mathscr{C}$. To this end, we arbitrarily choose a clause
$C_j\in\mathscr{C}$. Since the corresponding vertex $c_j$ is not
adjacent to any member of $\{s_2, s_4\}\cup\{v_i,v_i': 1\le i\le
n\}$, there exists some $i$ with $1\le i\le n$ such that $c_j$ is
dominated by $u_i\in D_t$ or $\bar{u}_i\in D_t$.

Suppose that $c_j$ is dominated by $u_i\in D_t$. Then $u_i$ is
adjacent to $c_j$ in $G$, that is, the literal $u_i$ is in the
clause $C_j$ by the construction of $G$. Since $u_i\in D_t$, we have
$t(u_i)=T$ by (\ref{e4.1}), which implies that $t$ satisfies the
clause $C_j$.

Suppose that $c_j$ is dominated by $\bar{u}_i\in D_t$. Then
$\bar{u}_i$ is adjacent to $c_j$ in $G$, that is, the literal
$\bar{u}_i$ is in the clause $C_j$. Since $\bar{u}_i\in D_t$, we
have $t(u_i)=F$ by (\ref{e4.1}), which implies that $\bar{u}_i$ is
assigned the truth value $T$ by $t$, so the clause $C_j$ is
satisfied by $t$.

The arbitrariness of $j$ with $1\le j\le m$ shows that all the
clauses in $\mathscr{C}$ is satisfied, that is, $\mathscr{C}$ is
satisfiable.

Conversely, suppose that $\mathscr{C}$ is satisfiable, and let $t:
U\to \{T,F\}$ be a satisfying truth assignment for $\mathscr{C}$.
Construct a subset $D'\subseteq V(G)$ as follows. If $t(u_i)=T$,
then put the vertex $u_i$ in $D'$; if $t(u_i)=F$, then put the
vertex $\bar{u_i}$ in $D'$. Clearly, $|D'|=n$. Since $t$ is a
satisfying truth assignment for $\mathscr{C}$, for each
$j=1,2,\ldots,m$, at least one of literals in $C_j$ is true under
the assignment $t$. It follows that the corresponding vertex $c_j$
in $G$ is adjacent to at least one vertex in $D'$ since $c_j$ is
adjacent to each literal in $C_j$ by the construction of $G$. Let
$D_t'=D'\cup \{s_2,s_4,v_1,\ldots,v_n\}$. Clearly, $D_t'$ is a
dominating set of $G$ and $|D_t'|=2n+2$. Since $s_2$ and $s_4$ are
dominated by each other, $u_i$ and $\bar{u}_i$ are dominated by
$v_i\in D_t'$ for each $i=1,2,\ldots,n$, $D_t'$ is also a total
dominating set of $G$. Hence, $\gamma_t(G)\leq |D_t'|=2n+2$. By
Claim 4.1, $\gamma(G)\geq 2n+2$. Therefore, $\gamma_t(G)=2n+2$.
\end{pf}

\item [Claim 4.3]
{\it For any $e\in E(G)$, $\gamma_t(G-e)\leq 2n+3$.}

\begin{pf}
We first assume $e=s_2s_3$ or $e=v_i\bar{u_i}$ for some $i$ with
$1\le i\le n$, and let $D_t'=(\cup_{i=1}^n\{u_i,v_i\})\cup
\{c_1,s_1,s_4\}$. It is easy to see that $D_t'$ is a total
dominating set of $G-e$. Secondly, assume $e=s_1c_j$ for some $j$
with $1\le j\le m$, and let $D_t'=(\cup_{i=1}^n\{u_i,v_i\})\cup
\{s_2,s_3,s_4\}$. Then $D_t'$ is a total dominating set of $G-e$.
Otherwise, let $D_t'=(\cup_{i=1}^n\{v_i,\bar{u_i}\})\cup
\{s_1,s_2,s_4\}$. Then $D_t'$ is a total dominating set of $G-e$.
Hence, $\gamma_t(G-e)\leq |D_t'|=2n+3$.
\end{pf}

\item [Claim 4.4]
{\it $\gamma_t(G)=2n+2$ if and only if $b_t(G)=1$.}

\begin{pf}
Assume $\gamma_t(G)=2n+2$ and take $e=s_2s_4$. Suppose that
$\gamma_t(G-e)=\gamma_t(G)$. Let $D_t'$ be a $\gamma_t$-set of
$G-e$. As $D_t'$ is also a $\gamma_t$-set of $G$, by Claim 4.1 we
have  $c_j\notin D_t'$ for every $j$ and $D_t'\cap
V(H)=\{s_2,s_4\}$, which contradicts the fact that $s_2$ and $s_4$
could not be dominated by each other in $G-e$. This contradiction
shows that $\gamma_t(G-e)>\gamma_t(G)$, whence $b_t(G)=1$.

Now, assume $b_t(G)=1$. By Claim 4.1, we have that $\gamma_t(G)\geq
2n+2$. Let $e'$ be an edge such that $\gamma_t(G-e')>\gamma_t(G)$.
By Claim 4.3, we have that $\gamma_t(G-e)\leq 2n+3$. Thus, $2n +
2\leq \gamma_t(G)< \gamma_t(G-e')\leq 2n+3$, which yields
$\gamma_t(G)=2n+2$.
\end{pf}

\end{description}

It follows from Claim 4.2 and Claim 4.4 that $b_t(G)=1$ if and only
if $\mathscr{C}$ is satisfiable. The theorem follows.
\end{pf}

\section{NP-hardness of reinforcement}

In this section, we will show that the problem determining the
reinforcements of general graphs is NP-hard. We first state it
as the following decision problem.

\begin{center}
\begin{minipage}{130mm}
\setlength{\baselineskip}{24pt}

\vskip6pt\noindent {\bf Reinforcement problem:}

\noindent {\bf Instance:}\ {\it A graph $G$ and a positive integer
$k$.}

\noindent {\bf Question:}\ {\it Is $r(G)\le k$?}

\end{minipage}
\end{center}

\vskip6pt\begin{thm} The reinforcement problem is NP-hard.
\end{thm}

\begin{pf}
The reinforcement problem is clearly in NP. In the following, we
show the NP-hardness of the reinforcement problem by reducing
the $3$-satisfiability problem to it in polynomial time.

Let $U=\{u_1,u_2,\ldots,u_n\}$ and $\mathscr{C}=\{C_1,C_2,
\ldots,C_m\}$ be an arbitrary instance of the $3$-satisfiability
problem. We will construct a graph $G$ and an integer $k$ such that
$\mathscr{C}$ is satisfiable if and only if $r(G)\leq k$. Such a
graph $G$ can be constructed as follows.

For each $i=1,2,\ldots,n$, corresponding to the variable $u_i\in U$,
associate a triangle $T_i$ with vertex-set $\{u_i,\bar{u_i},v_i\}$.
For each $j=1,2,\ldots,m$, corresponding to the clause
$C_j=\{x_j,y_j,z_j\}$, associate a single vertex $c_j$ and add edges
$(c_j,x_j), (c_j,y_j)$ and $(c_j,z_j)$, $1\le j\le m$. Finally, add
a vertex $s$ and join $s$ to every vertex $c_j$ and set $k=1$.

Figure~\ref{f3} shows an example of the graph obtained when
$U=\{u_1,u_2,u_3,u_4\}$ and $\mathscr{C}=\{C_1,C_2,C_3\}$, where
$C_1=\{u_1,u_2,\bar{u_3}\}, C_2=\{\bar{u_1},u_2,u_4\},
C_3=\{\bar{u_2},u_3,u_4\}$.

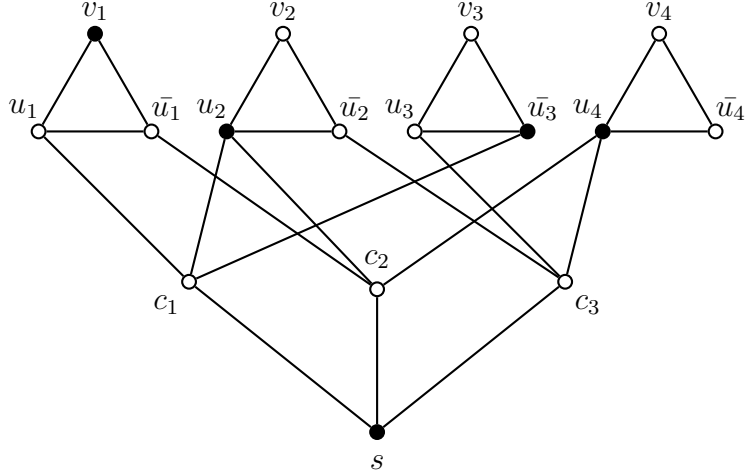
\begin{figure}[ht]
\begin{center}
\begin{pspicture}(-5,-.5)(5,6.7)

\cnode*(0,0){3pt}{s}\rput(0,-.4){$s$}

\cnode(0,1.9){3pt}{c2}\rput(0,2.3){$c_2$} \ncline{c2}{s}
\cnode(-2.5,2){3pt}{c1}\rput(-2.8,1.7){$c_1$} \ncline{c1}{s}
\cnode(2.5,2){3pt}{c3}\rput(2.8,1.7){$c_3$} \ncline{c3}{s}

\cnode(-4.5,4){3pt}{u1}\rput(-4.7,4.3){$u_1$}
\cnode(-3,4){3pt}{u1'}\rput(-2.8,4.3){$\bar{u_1}$} \ncline{u1}{u1'}
\cnode*(-2,4){3pt}{u2}\rput(-2.2,4.3){$u_2$}
\cnode(-0.5,4){3pt}{u2'}\rput(-0.3,4.3){$\bar{u_2}$} \ncline{u2}{u2'}
\cnode(0.5,4){3pt}{u3}\rput(0.3,4.3){$u_3$}
\cnode*(2,4){3pt}{u3'}\rput(2.2,4.3){$\bar{u_3}$} \ncline{u3}{u3'}
\cnode*(3,4){3pt}{u4}\rput(2.8,4.3){$u_4$}
\cnode(4.5,4){3pt}{u4'}\rput(4.7,4.3){$\bar{u_4}$} \ncline{u4}{u4'}
\cnode*(-3.75,5.3){3pt}{v1}\rput(-3.75,5.6){$v_1$} \ncline{v1}{u1} \ncline{v1}{u1'}
\cnode(-1.25,5.3){3pt}{v2}\rput(-1.25,5.6){$v_2$} \ncline{v2}{u2} \ncline{v2}{u2'}
\cnode(1.25,5.3){3pt}{v3}\rput(1.25,5.6){$v_3$} \ncline{v3}{u3} \ncline{v3}{u3'}
\cnode(3.75,5.3){3pt}{v4}\rput(3.75,5.6){$v_4$} \ncline{v4}{u4} \ncline{v4}{u4'}

\ncline{c1}{u1} \ncline{c1}{u2} \ncline{c1}{u3'}
\ncline{c2}{u1'} \ncline{c2}{u2} \ncline{c2}{u4}
\ncline{c3}{u2'} \ncline{c3}{u3} \ncline{c3}{u4}
\end{pspicture}
\caption{\label{f3}\footnotesize An instance of the reinforcement
problem resulting from an instance of the $3$-satisfiability
problem, in which $U=\{u_1,u_2,u_3,u_4\}$ and
$\mathscr{C}=\{\{u_1,u_2,\bar{u_3}\},\{\bar{u_1},u_2,u_4\},
\{\bar{u_2},u_3,u_4\}\}$. Here $k=1$ and $\gamma=5$, where the set
of bold points is a $\gamma$-set.}
\end{center}
\end{figure}

It is easy to see that the construction can be accomplished in
polynomial time. All that remains to be shown is that $\mathscr{C}$
is satisfiable if and only if $r(G)=1$. To this aim, we first prove
the following two claims.

\begin{description}

\item [Claim 5.1]
{\it $\gamma(G)=n+1$.}

\begin{pf}
Use the symbol $N[s]$ to denote the closed-neighborhood of $s$ in
$G$, that is, $N[s]=\{u\in V(G): us\in E\}\cup \{s\}$. On the one
hand, let $D$ be a $\gamma$-set of $G$, then $\gamma(G)=|D|\geq n+1$
since $|D\cap V(T_i)|\geq 1$ and $|D\cap N[s]|\geq 1$. On the other
hand, $D'=\{s,u_1,u_2,\ldots,u_n\}$ is a dominating set of $G$,
which implies that $\gamma(G)\le |D'|=n+1$. It follows that
$\gamma(G)=n+1$.
\end{pf}

\item [Claim 5.2]
{\it If there exists an edge $e\in E(\bar{G})$ such that
$\gamma(G+e)=n$, and let $D_e$ be a $\gamma$-set of $G+e$, then
$|D_e\cap V(T_i)|=1$ for each $i=1,2,\ldots,n$, while $c_j\notin
D_e$ for each $j=1,2,\ldots,m$.}

\begin{pf}
Suppose to the contrary that $|D_e\cap V(T_{i_0})|=0$ for some $i_0$
with $1\le i_0\le n$. Then one end-vertex of the edge $e$ should be
$v_{i_0}$ since $D_e$ dominates it via the edge $e$ in $G+e$, and
for every $i\neq i_0$, $|D_e\cap V(T_i)|\geq 1$ since $D_e$
dominates $v_i$. By the hypotheses, two literals $u_{i_0}$ and
$\bar{u}_{i_0}$ do not simultaneously appear in the same clause in
$\mathscr{C}$, they are not incident with the same vertex $c_j$ in
$G$ for some $j$. Since $u_{i_0}$ and $\bar{u}_{i_0}$ should be
dominated by $D_e$, there exist two distinct vertices $c_j, c_l\in
D_e$ such that $c_j$ dominates $u_{i_0}$ and $c_l$ dominates
$\bar{u}_{i_0}$. Thus, $|D_e|\geq n+1$, a contradiction. Hence,
$|D_e\cap V(T_i)|=1$ for each $i=1,2,\ldots,n$, and $c_j\notin D_e$
for every $j$ since $|D_e|=n$.
\end{pf}

\end{description}

We now show that $\mathscr{C}$ is satisfiable if and only if
$r(G)=1$.

Suppose that $\mathscr{C}$ is satisfiable, and let $t: U\to \{T,F\}$
be a satisfying truth assignment for $\mathscr{C}$. We construct a
subset $D'\subseteq V(G)$ as follows. If $t(u_i)=T$ then put the
vertex $u_i$ in $D'$; if $t(u_i)=F$ then put the vertex $\bar{u_i}$
in $D'$. Then $|D'|=n$. Since $t$ is a satisfying truth assignment
for $\mathscr{C}$, for each $j=1,2,\ldots,m$, at least one of
literals in $C_j$ is true under the assignment $t$. It follows that
the corresponding vertex $c_j$ in $G$ is adjacent to at least one
vertex in $D'$ since $c_j$ is adjacent to each literal in $C_j$ by
the construction of $G$. Without loss of generality let $t(u_1)=T$,
then $D'$ is a dominating set of $G+su_1$, and hence
$\gamma(G+su_1)\leq |D'|=n$. By Claim 5.1, we have $\gamma(G)=n+1$.
It follows that $\gamma(G+su_1)\leq n<n+1=\gamma(G)$, which implies
$r(G)=1$.

Conversely, assume $r(G)=1$. Then there exists an edge $e$ in
$\bar{G}$ such that $\gamma(G+e)=n$. Let $D_e$ be a $\gamma$-set of
$G+e$. By Claim 5.2, $|D_e\cap V(T_i)|=1$ for each $i=1,2,\ldots,n$, 
and $c_j\notin D_e$ for each $j=1,2,\ldots,m$. Define $t: U\to \{T,F\}$ by
 \begin{equation}\label{e5.1}
 t(u_i)=\left\{
\begin{array}{ll}
 T \ & {\rm if}\ u_i\in D_e \ {\rm or} \ v_i\in D_e, \\
 F \ & {\rm if}\ \bar {u_i}\in D_e,
\end{array}
 \right.
 \ i=1,2,\ldots,n.
 \end{equation}

We will show that $t$ is a satisfying truth assignment for
$\mathscr{C}$. It is sufficient to show that every clause in
$\mathscr{C}$ is satisfied by $t$. 

Consider arbitrary clause $C_j\in\mathscr{C}$ with $1\le j\le m$. By
Claim 5.2, the corresponding vertex $c_j$ in $G$ is dominated by
$u_i$ or $\bar{u}_i$ in $D_e$ for some $i$. Suppose that $c_j$ is
dominated by $u_i\in D_e$. Then $u_i$ is adjacent to $c_j$ in $G$,
that is, the literal $u_i$ is in the clause $C_j$ by the
construction of $G$. Since $u_i\in D_e$, we have $t(u_i)=T$ by
(\ref{e5.1}), which implies that $C_j$ is satisfied by $t$. Suppose
that $c_j$ is dominated by $\bar{u}_i\in D_e$. Then $\bar{u}_i$ is
adjacent to $c_j$ in $G$, that is, the literal $\bar{u}_i$ is in the
clause $C_j$. Since $\bar{u}_i\in D_e$, we have $t(u_i)=F$ by
(\ref{e5.1}), which implies that $\bar{u}_i$ is assigned the truth
value $T$ by $t$, so the clause $C_j$ is satisfied. The
arbitrariness of $j$ with $1\le j\le m$ shows that all the clauses
in $\mathscr{C}$ is satisfied by $t$, that is, $\mathscr{C}$ is
satisfiable.

The theorem follows.
\end{pf}

\end{document}